\documentclass[11pt]{article}

\usepackage{amsmath,latexsym,epsfig}
\usepackage{amsfonts,amssymb,amscd}

\title{\bf The number of ramified coverings of the sphere by the torus  and surfaces of higher 
genera}

\author{
I.P.Goulden\thanks{Dept. of Combinatorics and Optimization,
University of Waterloo, Waterloo, Ontario, Canada},
D.M.Jackson\thanks{Dept. of Combinatorics and Optimization,
University of Waterloo, Waterloo, Ontario, Canada} and
A.Vainshtein\thanks{Dept. of Mathematics and Dept. of Computer Science, University of Haifa, Haifa, Israel}
}

\date{January 11, 1999}

\begin{document}
 \maketitle

 \newtheorem{theorem}{Theorem}[section]
 \newtheorem{proposition}[theorem]{Proposition}
 \newtheorem{definition}[theorem]{Definition}
 \newtheorem{axiom}[theorem]{Axiom}
 \newtheorem{lemma}[theorem]{Lemma}
 \newtheorem{corollary}[theorem]{Corollary}
 \newtheorem{remark}[theorem]{Remark}
 \newtheorem{example}[theorem]{Example}
 \newtheorem{conjecture}[theorem]{Conjecture}

\def\cA{{\cal{A}}}
\def\cB{{\cal{B}}}
\def\cC{{\cal{C}}}
\def\cR{{\cal{R}}}
\def\cS{{\cal{S}}}
\def\cT{{\cal{T}}}

\def\bfp{{\rm\bf p}}
\def\bfw{{\rm\bf w}}
\def\bfx{{\rm\bf x}}

\def\sM{{\mathbb M}}
\def\sP{{\sf P}}
\def\sS{{\mathbb S}}

\def\symgp{{\mathfrak S}}
\def\rats{{\mathbb Q}}

\def\proof{{\rm\bf Proof:\quad}}
\def\qed{{\hfill{\large $\Box$}}}

\def\mybar{\overline}
\def\sym#1{\Theta_{#1}}

\def\atp#1#2{\stackrel{\scriptstyle{#1}}{\scriptstyle{#2}}}

\def\aatp#1#2#3{\atp{\atp{#1}{#2}}{#3}}


\begin{abstract}
We obtain an explicit expression for the number of ramified coverings of the sphere
by the torus with given ramification type for a small number of ramification points,
and conjecture this to be true for an arbitrary  number of ramification points.
In addition, the conjecture is proved for simple coverings of the sphere by the torus.
We obtain corresponding expressions for surfaces of higher genera
for small number of ramification points, and conjecture the general form for this
number in terms of a symmetric polynomial that appears to be new.
The approach involves the analysis of the action of a transposition to derive
a system of linear partial differential equations that give the generating series for
the desired numbers.  
\end{abstract}

\section{Introduction}\label{Si}


\subsection{Background and new results}\label{ssback}
Let $f\colon\sM\rightarrow\sS^2$ be a non-constant meromorphic function on a compact connected Riemann
surface $\sM$ of genus $g\ge0.$ Then there exists an integer $n\ge1$, called the {\em degree} of $f,$
such that $\vert f^{-1}(p)\vert=n$ for all but a finite number of points $p\in\sS^2$ called {\em critical values};
$f$ is called a {\em ramified $n$-fold covering} of $\sS^2.$

There are two definitions of topological equivalence of ramified coverings. In one of them, two
coverings $f_1$ and $f_2$ are considered equivalent if there exist homeomorphisms
$\pi\colon\sM\rightarrow\sM$ and $\rho\colon\sS^2\rightarrow\sS^2$ such that $\rho f_1=f_2\pi.$
The equivalence classes under this definition correspond to the connected components of the space of meromorphic
functions. For the description and enumeration of such equivalence classes see~\cite{KZ,N}

In the second, more restrictive definition, two coverings $f_1$ and $f_2$ are considered equivalent if there 
exists a homeomorphism $\pi\colon\sM\rightarrow\sM$ such that $f_1=f_2\pi.$
In this paper we address the problem of determining the number of equivalence classes of ramified coverings 
under this more restrictive equivalence. This problem arose originally in the classic paper of Hurwitz~\cite{Hrfgv},
and is often called {\em Hurwitz's Enumeration Problem}.

Recall that a critical value $p$ of $f$ is called {\em simple} if  $\vert f^{-1}(p)\vert=n-1.$
A covering $f$ is called {\em simple} if all of its critical values are simple, and {\em almost simple}
if at most one critical value, which is assumed to be equal to $\infty,$ is non-simple.
The preimages of $\infty$ are called the {\em poles} of $f.$
Let $\alpha_1\ge\ldots\ge \alpha_m\ge1$ be the orders of the poles;
since $\alpha_1+\cdots+\alpha_m=n,$ we get a partition $\alpha$
of $n,$ which is indicated by $\alpha\vdash n.$ The number of parts is denoted by $\l(\alpha).$
The specific question addressed in this paper is that of giving an explicit expression for the
number $\mu_m^{(g)}(\alpha)$ of almost simple ramified $n$-fold coverings with a prescribed distribution $\alpha$ of
ramification orders. We call $\alpha$ the {\em ramification type} of $f.$

A general answer to this question been given by a number of authors
(see for example~\cite{Hrfgv,medn}).
These answers allow us to determine $\mu_m^{(g)}(\alpha)$ in principle, but do not
give explicit information, and are largely intractable because of the character sums they contain.
For example, the following generating series~(see~\cite{GJtransf}), in principle, gives complete
information about the $\mu_m^{(g)}(\alpha).$
It gives $\mu_m^{(g)}(\alpha)/(n+m+2g-2)!$ as the coefficient of $z^g u^{n+m+2g-2}p_\alpha$ in
\begin{eqnarray*}
\log\left( \sum_{n\ge0}\sum_{\theta\vdash n}\frac{1}{n!} \chi^\theta_{[1^n]}
s_\theta e^{u\eta(\theta)}\right).
\end{eqnarray*}
In this series $p_k$ is the power sum symmetric function of degree $k,$ $p_\alpha=p_{\alpha_1}p_{\alpha_2}\cdots,$
$\chi^\theta$ is the character of the ordinary irreducible representation of $\symgp_n$ 
(the symmetric group on $n$ symbols) indexed by $\theta,$
$$s_\theta =\frac{1}{n!}\sum_{\alpha\vdash n}\vert\cC_\alpha\vert  \chi^\theta_\alpha p_\alpha,$$
the Schur function indexed by $\theta,$ $\vert\cC_\alpha \vert$ is the size of the conjugacy
class of $\symgp_n$  indexed by $\alpha,$
and
$$\eta(\theta)=\sum_i \binom{\theta_i}{2} - \sum_i \binom{\tilde{\theta_i}}{2},$$
where $\theta=(\theta_1,\ldots),$ $\tilde{\theta}=(\tilde{\theta}_1,\ldots),$ 
the conjugate of $\theta.$

A different approach to this problem comes from singularity theory and is based on the notion of
the Lyashko-Looijenga map (see~\cite{A1}). This map takes a meromorphic function $f$
to the polynomial whose roots are the critical values of $f.$ It was shown in~\cite{L} that
for $\sM=\sS^2$ and the generic polynomial, $f,$ the Lyashko-Looijenga map is  a finite covering
whose degree coincides, up to a factor of $n,$ with $\mu_1^{(0)}((n)).$
This approach was extended to generic Laurent polynomials on the sphere in~\cite{Arntcc},
and to generic meromorphic functions on the sphere in~\cite{GL}.
The case of polynomials on $\sS^2$ with several non-simple critical values is treated in~\cite{Z}.

A further approach to the problem comes from enumerative algebraic geometry. Recall that the 
characteristic number $R_d(a,b)$ (with $a+b=3d-1$) is defined as the number of irreducible degree $d$
rational curves passing through $a$ fixed general points, and tangent to $b$ fixed general lines.
Kontsevich and Manin~\cite{K,KM} were the first to relate the characteristic numbers
(in the case $b=0$) to Chern classes of the moduli space of stable maps and to give recursions for
these numbers. Their results were later extended in~\cite{P} to the case of arbitrary $b.$
A further extension to the characteristic numbers $R_d(a,b,c)$ of curves that satisfy
the above conditions and, in addition, are tangent to $c$ fixed general lines at
fixed general points of tangency can be found in~\cite{EK}.
The same approach allows the finding of recurrence equations of the number $\mu^{(g)}_n(1^n)$
of simple coverings, at least in the cases $g=0$ and $g=1$ (see~\cite{V}).
The case of general almost simple coverings has not yet been studied from this point of view.
The first steps in this direction have been carried out in~\cite{ELSV}, where $\mu^{(g)}_n(\alpha)$
is related
to the top Segre class of  a certain bundle over the moduli space of stable curves.

It is convenient to rescale  $\mu_m^{(g)}(\alpha)$ by writing it in terms of a quantity,  $f_m^{(g)}(\alpha),$
which we  introduce in this paper, where
\begin{eqnarray}\label{newequ}
\mu_m^{(g)}(\alpha)=\frac{1}{n!}\vert\cC_\alpha\vert\,(n+m+2g-2)! \prod_{j=1}^m
\frac{\alpha_j^{\alpha_j}}{(\alpha_j-1)!}f_m^{(g)}(\alpha).
\end{eqnarray}
Much of this paper is concerned with determining explicit expressions and properties
of $f_m^{(g)}(\alpha)$ defined by this rescaling.
In terms of  this rescaling, the present state of knowledge of $\mu_m^{(g)}(\alpha)$ is summarized as follows. 
\begin{eqnarray}
f_m^{(0)}(\alpha) &=& n^{m-3}, \quad\alpha\vdash n, l(\alpha)=m,  \label{er1} \\
f_1^{(g)}((n)) &=& \frac{1}{2^{2g}}n^{2g-2}
\left[x^{2g}\right]\,\left(\frac{\sinh x}{x}\right)^{n-1},\quad n\ge1,  \label{er2} \\
f_2^{(1)}((n-r,r)) &=& \frac{1}{24}
\left(n^2-(r+1)n+r^2\right), \quad r=1,2,3, \,\, n-r\neq r,\label{er3}
\end{eqnarray}
where $[x^r]g(x)$ is the coefficient of $x^r$ in $g(x).$ 

Although it has  a strikingly simple form, result~(\ref{er1}), which is due to 
 Hurwitz~\cite{Hrfgv} (see also Goulden and Jackson~\cite{GJtransf},
Strehl~\cite{vs} and Goryunov and Lando~\cite{GL}), is a remarkable one.
The case $\alpha=[1^n]$ was  rediscovered recently by Crescimanno and Taylor~\cite{CT}.
Result~(\ref{er2}) was obtained by Shapiro, Shapiro and Vainshtein~\cite{ssv},
who also showed that it could be obtained from a result of Jackson~\cite{dmjcppcc} on
ordered factorizations of permutations into transpositions.
Result~(\ref{er3}) was obtained by Shapiro, Shapiro and Vainshtein~\cite{ssv}.

Thus, for the {\em sphere}, complete information is known; 
for the {\em torus} information is known only for ramification types $(n)$ and $(n-r,r),$ $r=1,2,3$ and
$n-r>r$;  and for {\em all surfaces},  explicit information is known only for one ramification point.

In this paper we determine several  new explicit results for $\mu_m^{(g)}(\alpha).$ 
First, we find the number of simple coverings of a sphere by a torus, thus generalizing
Hurwitz's result for the sphere.
\begin{theorem}\label{T11}
For $n\ge1,$
$$f^{(1)}_n(1^n)=\frac{1}{24}\left(n^n-n^{n-1}-\sum_{i=2}^n\binom{n}{i}(i-2)!n^{n-i}\right).$$
\end{theorem}

Next, we consider almost simple coverings of the sphere by a torus. We prove that result~(\ref{er3})
remains valid for arbitrary coverings with two ramification points.
\begin{theorem}\label{T12}
For $n\ge1$ and $0<r<n,$
$$f_2^{(1)}((n-r,r)) = \frac{1}{24}
\left(n^2-(r+1)n+r^2\right).$$
\end{theorem}

Moreover, we find explicit expressions for the number of almost simple coverings
of a sphere by a torus with up to 6 ramification points (see Appendix).
These expressions, together with Theorem~\ref{T11}, allow us to conjecture the general explicit result
that gives complete information for the torus, and that can therefore be regarded as an extension
of result~(\ref{er1}).

\begin{conjecture}\label{cfm24} For $m\ge 1,$
\begin{eqnarray}\label{efm24}
f^{(1)}_m(\alpha)=\frac{1}{24}\left(n^m-n^{m-1}-\sum_{i=2}^m(i-2)!e_{i}n^{m-i}\right),
\end{eqnarray}
where $e_i$ is the $i$-th elementary symmetric function in $\alpha_1,\ldots,\alpha_m,$
and $e_1=n=\alpha_1+\cdots+\alpha_m.$
\end{conjecture}

Furthermore, we determine new explicit results for $\mu^{(g)}_m(\alpha)$ in the cases
$m+g\le 6,$ $m,g\ge2.$ These results consist of explicit expressions for the corresponding
$f^{(g)}_m(\alpha)$ and are collected together in the Appendix. Each result $f^{(g)}_m(\alpha)$ is 
a symmetric polynomial in $\alpha_1,\ldots,\alpha_m.$ We have found that these polynomials are
expressed more compactly in terms of the elementary symmetric functions $e_k,$ $k=1,\ldots,m,$
and, for this reason, this is the presentation of $f^{(g)}_m(\alpha)$ that is given
in the Appendix ({\em cf.} also Conjecture~\ref{cfm24}).
In all cases we also observe that  the total degree of $f_m^{(g)}(\alpha)$ 
is $m+3g-3$ and we therefore make the following conjecture.
\begin{conjecture}\label{conject1}
For $g=0, m\ge3,$ and $m\ge1, g\ge1,$
$f_m^{(g)}(\alpha_1,\ldots,\alpha_m)$ is a symmetric polynomial in $\alpha_1,\ldots,\alpha_m$
of total degree $m+3g-3.$
\end{conjecture}

Conjecture~\ref{conject1} is in agreement also with the previously known results
(\ref{er1}), (\ref{er2}), and (\ref{er3}). Note that for $\alpha\vdash n,$ we have $n=e_1,$ so 
Hurwitz's result~(\ref{er1}) is rewritten in these
terms as $f_m^{(0)}=e_1^{m-3},$ which is a (symmetric) polynomial for $m\ge3.$
For result~(\ref{er2}), we have
\begin{eqnarray*}
\left[x^{2g}\right] \left(\frac{\sinh x}{x}\right)^{n-1} 
&=& \left[t^{g }\right]\,\left( \sum_{i\ge0}(-1)^i\frac{t^i}{(2i+1)!}\right)^{n-1} \\
&=& \sum_{\atp{i_1,\ldots,i_g\ge0}{i_1+2i_2+\cdots+gi_g=g}} (n-1)_{i_1+\cdots+i_g}\prod_{j=1}^g
\frac{1}{i_j!(2j+1)!^{i_j}},
\end{eqnarray*}
where $(u)_i=u(u-1)\cdots(u-i+1)$ for a nonnegative integer $i.$
Thus $f^{(g)}_1$ is a polynomial of degree $3g-2$ in $n=e_1,$ 
in agreement with Conjecture~\ref{conject1} for $g\ge1.$

Although we have been unable to use this approach to prove
Conjecture~\ref{cfm24} and Conjecture~\ref{conject1}, 
we believe that the symmetric functions  $f_m^{(g)}$ are of considerable interest and 
would reward further study. Explicit expressions for a few of them are given in the Appendix in 
terms of the elementary symmetric functions. It may be that there is a more natural basis
for the   $f_m^{(g)},$ but we have been unable to determine it.


\subsection{New methods and the organization of the paper}
Our approach uses Hurwitz reduction of the question to the ordered factorization of permutations in $\symgp_n$
into transpositions that generate $\symgp_n.$ 
We obtain a differential equation for the generating series by an analysis that
lies at the centre of  Hurwitz's approach. Our point of departure from Hurwitz's approach 
is to show that the symmetrization
of this series satisfies a differential equation in a new set of variables, and that its solution
is a rational function in the transformed variables. This enables us, not only  to obtain
additional explicit results that extend what is currently known, but also to make conjectures about
the general form of $\mu_m^{(g)}(\alpha).$ 

In Section~\ref{stof} the determination of $\mu_m^{(g)}(\alpha)$ is expressed in terms of 
ordered factorizations of permutations in $\symgp_n$ into transpositions
such that the transpositions  generate
$\symgp_n.$ A differential equation for the generating series for the number $c_g(\alpha)$
of such factorizations is  given in Section~\ref{scaja} by analyzing the action of
transpositions combinatorially. Section~\ref{ssfP} gives a partial differential
equation induced by symmetrizing the portion of the generating series relating to given genus $g$ and 
a prescribed number $m$ of ramification points. 
A general form for the solution of this system is
conjectured in Section~\ref{soi}. 
Section~\ref{scomments} gives brief details about the computation of the results in the Appendix,
and a proof of Theorem~\ref{T12}.
Section~\ref{ST11} gives a proof of Theorem~\ref{T11} and makes use of one of the new variables
that was introduced in Section~\ref{ssfP}.



\section{Transitive ordered factorizations}\label{stof}
The question
of determining the number $\mu_m^{(g)}(\alpha)$ of almost simple ramified $n$-fold coverings,
by a surface of genus $g,$ with a prescribed  ramification type $\alpha$,
can be translated into one concerning products of transpositions, and the construction 
is given by Hurwitz~\cite{Hrfgv}.
Let $c_g(\alpha)$ be the number of $j$-tuples of permutations $(\sigma_1,\ldots,\sigma_j)$
such that, for an arbitrary but fixed $\pi\in\cC_\alpha,$ the conjugacy class of $\symgp_n$ indexed
by $\alpha\vdash n,$ the following conditions hold:
\begin{enumerate}
   \item $\pi=\sigma_1\cdots\sigma_j,$ \label{d1}
   \item $\sigma_1,\ldots,\sigma_j\in\cC_{[2\,1^{n-2}]},$  \label{d2}
   \item  $\sigma_1,\ldots,\sigma_j$ generate $\symgp_n,$  \label{d3}
   \item $j=\mu(\alpha)+2g,$ where $\mu(\alpha)=n+l(\alpha)-2.$  \label{d4}
\end{enumerate}
We note that $\mu(\alpha)$ is the minimum value of $j$ for which conditions (\ref{d1}), (\ref{d2}) 
and~(\ref{d3}) hold. We call $(\sigma_1,\ldots,\sigma_r)$ an {\em ordered transitive factorization}
of $\pi$ into transpositions, and is said to be {\em minimal} if $g=0.$
The numbers $\mu_m^{(g)}(\alpha)$ and $c_g(\alpha)$ are related by

\begin{eqnarray}\label{emug}
\mu_m^{(g)}(\alpha)  =\frac{1}{n!}\vert\cC_\alpha\vert c_g(\alpha).
\end{eqnarray}
Condition~\ref{d1} ensures that the ramification type is $\alpha=[1^{a_1},2^{a_2},\ldots],$
where $a_k$ is the number of ramification points of order $k.$
Condition~\ref{d2} ensures that the ramification points are simple.
Condition~\ref{d3} ensures that the surface is connected. 
Condition~\ref{d4} ensures that the sphere is covered by a surface of genus $g.$

It is this reformulation of the original question that we now consider in detail.
Goulden~\cite{Gdosf} showed that such  problems can be considered by 
analyzing the combinatorial action 
of a transposition on a permutation through the use of differential operators
acting on monomials that encode the cycle structure of a permutation.
These techniques were 
extended in \cite{GJtfst} to factorizations into cycles of length $k$ alone.
We refer to this technique as the {\em cut-and-join} analysis of the action of a $k$-cycle.
The cut-and-join analysis is applied to a graphical encoding of transitive ordered factorizations
as connected edge-labelled and vertex-labelled graphs. These  appear in the work of Arnold~\cite{Arntcc},
and independently in~\cite{Gdosf,GJtfst}, and are called the {\em monodromy graph} of the ordered 
factorization.
For the question addressed in this paper it is necessary to relax the condition of minimality
in the cut-and-join analysis.

\section{A partial differential equation}\label{scaja}
Let $p_1,p_2,\ldots$ be indeterminates and  $\bfp=(p_1,p_2,\ldots).$
For $\alpha=(\alpha_1,\ldots,\alpha_m)$ where $\alpha_1,\ldots,\alpha_m$ are positive integers, 
let $p_\alpha=p_{\alpha_1}\cdots p_{\alpha_m}.$ Now consider the  generating series
\begin{eqnarray}\label{dPhi}
\Phi(u,z,\bfp)=
\sum_{\atp{n,m\ge1} {g\ge0}}
\sum_{\atp{\alpha\vdash n} {l(\alpha)=m}}
\frac{\vert\cC_\alpha\vert}{n!} c_g(\alpha) \frac{u^{n+m+2(g-1)}}{(n+m+2(g-1))!}z^g p_\alpha
\end{eqnarray}
for ordered transitive factorizations.
In Lemma~\ref{LPhi} we show that
 $\Phi$ satisfies  a second order partial differential equation of the second degree.

The core of the derivation is a cut-and-join analysis. Suppose that $\sigma\in\symgp_n$  is
the transposition that interchanges $s$ and $t,$ and that $\pi\in\symgp_n.$
Then there are two cases for the action of $\sigma$ on $\pi$ in the product $\pi\sigma\colon$
i)~if $s,t$ are on the same cycle of $\pi,$ then the cycle is cut into two cycles in $\pi\sigma,$
with $s,t$ on different cycles (here $\sigma$ is called a {\em cut} for $\pi$);
ii)~if $s,t$ are on different cycles of $\pi,$ then these cycles are joined to form a single cycle
in $\pi\sigma,$ which contains both $s$ and $t$ (here $\sigma$ is called a {\em join} for $\pi$).
This type of  analysis has been previously used in~\cite{GJtransf} and~\cite{GJtfst}, to obtain partial
differential equations for  other generating series   for ordered factorizations into transpositions.
In the proof of Lemma~\ref{LPhi} some of the details are suppressed because of the similarity
to this previous work.

\begin{lemma}\label{LPhi}
\begin{eqnarray}\label{eLPhi}
\frac{\partial\Phi}{\partial u} =
\frac{1}{2}\sum_{i,j\ge1}\left(
 ij p_{i+j}z \frac{\partial^2\Phi}{\partial p_i\partial p_j}
+  ij \,p_{i+j}\frac{\partial\Phi}{\partial p_i}\frac{\partial\Phi}{\partial p_j}
+(i+j)p_ip_j\frac{\partial\Phi}{\partial p_{i+j}}\right).
\end{eqnarray}
\end{lemma}
\proof  
Consider an ordered  transitive factorization $(\sigma_1,\ldots,\sigma_r)$ of $\pi$ that satisfies 
conditions~\ref{d1},~\ref{d2},~\ref{d3} and ~\ref{d4}, for any $n\ge1.$
Remove $\sigma_r$  but leave $\pi$ unchanged.  This is an ordered transitive factorization of $\pi$ with
its rightmost factor deleted.
Then the generating series for the set of such objects is
$$\frac{\partial\Phi}{\partial u},$$
since $\Phi$ is an exponential generating series in $u.$
This gives the left hand side of the partial differential equation.

To obtain the right hand side of the partial differential equation, we now
show that ordered transitive factorizations of $\pi$ with
the rightmost factor deleted
can be determined in another way. Let $G=G_{(\sigma_1,\ldots,\sigma_r)}$ be
the monodromy graph of $(\sigma_1,\ldots,\sigma_r)$: 
 the graph with vertex 
labels $\{1,\ldots,n\},$ and edge labels $\{1,\ldots,r\},$ in which the edge labelled $i$ joins
vertices interchanged by $\sigma_i.$ Then condition~\ref{d3} implies that $G$ is a 
connected graph. 
 Let $T$ be the spanning tree of $G$ constructed by considering the edges in increasing order of their labels,
and selecting the edge labelled $i$ if and only if its incident vertices are in different components
of $G_{(\sigma_1,\ldots,\sigma_{i-1})}$ (this is Kruskal's Algorithm).  Then, by construction,
the transpositions corresponding to the edges of $T$ are joins.
Let $\sigma_r=(s,t).$ There are three cases. 

i) $\sigma_r$ corresponds to an edge of $G$ not in $T,$ and $\sigma_r$ is a cut for $\pi\sigma_r.$
Thus $(\sigma_1,\ldots,\sigma_{r-1})$ is an ordered transitive factorization of  $\pi\sigma_r,$
with $s$ and $t$ on the same cycle of $\pi\sigma_r.$
The contribution from this case is
  $$\frac{1}{2}\sum_{i,j\ge1} (i+j)p_ip_j\frac{\partial\Phi}{\partial p_{i+j}},$$
since the effect of multiplying by $\sigma_r$ is to cut a $i+j$-cycle into
an $i$-cycle and a $j$-cycle.

ii)  $\sigma_r$ corresponds to an edge of $G$ not in $T,$ and $\sigma_r$  is a join for $\pi{\sigma_r}.$.
Thus $G_{(\sigma_1,\ldots,\sigma_{r-1})}$ is an ordered transitive factorization 
of $\pi\sigma_r$ with $s$ and $t$ on different cycles of $\pi\sigma_r.$ The
 contribution from this case is 
$$\frac{1}{2}\sum_{i,j\ge1} ij p_{i+j}z \frac{\partial^2\Phi}{\partial p_i\partial p_j},$$
since the effect of multiplying by $\sigma_r$ is to join an $i$-cycle and a $j$-cycle,
to create an $i+j$-cycle.

iii) $\sigma_r$ corresponds to an edge of $T$; then $\sigma_r$  corresponds  a join for
$\pi\sigma_r$ and $G_{(\sigma_1,\ldots,\sigma_{r-1})}$ has exactly two components,
with $s$ and $t$ on different components. The contribution from this case is
$$\frac{1}{2}\sum_{i,j\ge1} ij p_{i+j}\frac{\partial\Phi}{\partial p_i}\frac{\partial\Phi}{\partial p_j},$$
where here the $i$-cycle and $j$-cycle that have been joined   come from each of two ordered
transitive factorizations (one for each component).
Note that $\Phi$ is an exponential generating series in both $u$ and vertices,
the latter being indicated by the division by $n!$ in the definition of $\Phi,$
so the product
$ij\frac{\partial\Phi}{\partial p_i}\frac{\partial\Phi}{\partial p_j}$
gives the correct cardinalities through the product of exponential generating series.
The result follows by combining these cases and equating the two expressions.
\qed

In~\cite{GJtransf} the special case of $\Phi$ corresponding to minimal ordered transitive
factorizations was considered. This is obtained by considering only $g=0,$ or, equivalently, by
setting $z=0$ in $\Phi.$ The partial differential equation obtained by setting $z=0$ in~(\ref{eLPhi})
was derived in~\cite{GJtransf} by a cut-and-join analysis, and an explicit solution was obtained. Although we have been
unable to obtain explicit solutions to~(\ref{eLPhi}) itself, we have been
able to obtain results for various positive values of $g$ by considering a symmetrized form of~(\ref{eLPhi}),
and the rest of the paper is devoted to this analysis.

\def\mv{\qquad}
\def\mvv{\qquad\qquad\qquad}

\section{The symmetrical form of $\Phi$}\label{ssfP}

We begin by defining a symmetrical form of $\Phi.$  For this purpose, let   
$\bfx^\alpha=x_1^{\alpha_1}\cdots x_m^{\alpha_m}$ and let
$\omega\alpha=(\alpha_{\omega(1)},\ldots,\alpha_{\omega(m)})$ where $\omega\in\symgp_m.$
Let 
\begin{eqnarray*} \Xi^{\{1,\ldots,m\}}p_\alpha = \left\{
\begin{array}{ccccc} 
\sum_{\pi\in\symgp_m}\bfx^{\pi(\alpha)} &\mbox{ if } l(\alpha)=m,\\
0&\mbox{otherwise,}
\end{array} \right.
\end{eqnarray*}
extended linearly to all series in the $p_i$'s.
Let
\begin{eqnarray}\label{eXi}
\Psi_m^{(g)}(\bfx)=\left.[z^g]\Xi^{\{1,\ldots,m\}} \Phi(u,z,\bfp)\right\vert_{u=1}.
\end{eqnarray}

Thus $\Psi_m^{(g)}(\bfx)$ is obtained by taking the term of degree $g$ in $z$ and total degree $m$ in the $p_i$'s
from $\Phi,$ then setting $u=1$ and symmetrizing the monomials in the $p_i$'s, replacing the subscripts of $p_i$'s
by exponents of $x$'s when we do so.
Of course, $c_g(\alpha)$ is recoverable from this series in a straightforward way.

In this section we determine a partial differential equation for $\Psi_m^{(g)}(\bfx)$ that is 
induced by the partial differential equation for $\Phi,$ given in  Lemma~\ref{LPhi}.
It is convenient to
express the equations in terms of $x_1,\ldots,x_m$ and, in addition, $w_1,\ldots,w_m,$ where
\begin{eqnarray}\label{wxe}
w_i=x_ie^{w_i}
\end{eqnarray}
has a unique solution $w_i\equiv w_i(x_i)$ as a power series in $x_i.$ Then
$w_1,\ldots,x_1,\ldots$ are algebraically independent, 
and it is clear, by differentiating the functional equation~(\ref{wxe}), that 
\begin{eqnarray}\label{yxw}
x_i\frac{\partial}{\partial x_i} = \frac{w_i}{1-w_i} \frac{\partial}{\partial w_i}.
\end{eqnarray}
The reason for introducing the $w_i$'s is that in~\cite{GJtransf} and~\cite{GJtfst}  
we showed that Hurwitz's result~(\ref{er1}) for $g=0$ can be expressed as
\begin{eqnarray}\label{Pe1m}
\Psi_m^{(0)}(\bfx) =\left( x_1\frac{\partial}{\partial x_1}+\cdots+x_m\frac{\partial}{\partial x_m}
\right)^{m-3}V_m(\bfw)
\end{eqnarray}
in terms of the new indeterminates, where
\begin{eqnarray}\label{Vkpw}
V_m(\bfw)=\prod_{i=1}^m \frac{w_i}{1-w_i} 
\end{eqnarray}
and $\bfw=(w_1,\ldots,w_m).$ This expression will be used as an initial condition in 
the partial differential equation for $\Psi_m^{(g)}(\bfx).$

The following mapping is needed  in the statement of the differential equation,
for expressing the action of $\Xi^{\{1,\ldots,m\}}$ in terms involving
the $p_i$'s. Let $f$ be a series in $x_1,\ldots,x_{m},$ and let $0\le i\le m-1.$ Then the mapping
$\sym{i}$ is defined by
$$\sym{i} f(x_1,\ldots,x_{m})=\sum_{\cR,\cS,\cT}f(\bfx_\cR,\bfx_\cS,\bfx_\cT)$$
where the sum is over all ordered partitions $(\cR,\cS,\cT)$ of $\{1,\ldots,m\}$
with $\vert\cR\vert=1,$  $\vert\cS\vert=i,$  $\vert\cT\vert=m-i-1,$ and where
$$(\bfx_\cR,\bfx_\cS,\bfx_\cT)=(x_{r_1},x_{s_1},\ldots,x_{s_i},x_{t_1},\ldots,x_{t_{m-i-1}}),$$
in which $s_1<\ldots<s_i,$ and $t_1<\ldots<t_{m-i-1}.$

The next three lemmas are quite technical, and
will be used for determining the action of $\Xi^{\{1,\ldots,m\}}$ on products
of the $p_i$'s. For this purpose, let
$$\Xi^{\{a_1,\ldots,a_m\}}p_\alpha=\left. \Xi^{\{1,\ldots,m\}}p_\alpha\right\vert_{x_i\mapsto x_{a_i}, i=1,\ldots,m},$$
where $a_1<\ldots<a_m.$

\begin{lemma}\label{L1}
Let  $\alpha,\beta$ be partitions with $l(\alpha)=k$ and  $l(\beta)=m.$ Then
$$\Xi^{\{1,\ldots,m+k\}}p_\alpha  p_\beta
=\sum_{(\cA,\cB)} \left( \Xi^{\cA} p_\alpha\right)  \left( \Xi^{\cB}p_\beta\right),$$
where the sum is over all ordered partitions $(\cA,\cB)$ of $\{1,\ldots,m+k\}$ with
$\vert\cA\vert=k$ and $\vert\cB\vert=m.$
\end{lemma}
\proof Immediate. \qed

\begin{lemma}\label{L2}
Let $\alpha$ be a partition, let $l(\alpha)=m,$ and let $1\le l\le m.$ Then
$$ \sum_{i\ge1} x_l^i \Xi^{\{1,\ldots,m\} - \{l\}} i\frac{\partial  p_\alpha}{\partial p_i}=
x_l\frac{\partial }{\partial x_l}\,\Xi^{\{1,\ldots,m\}}p_\alpha.$$
\end{lemma}
\proof Let $\alpha=(\alpha_1,\ldots,\alpha_m).$ It is convenient to regard $p_\alpha$ as the word
$p_{\alpha_1}\cdots p_{\alpha_m}$ in which each symbol is written with exponent equal to one.
Now apply the operator $\sum_{i\ge1}x_l^i i\partial/ \partial p_i$  to this product of $m$ $p_i$'s. Its 
effect, by the product rule, is to create a sum of $m$ terms, in which each of the $m$ $p_j$'s in the word
$p_\alpha$ is replaced, in turn, by $j x_l^j.$ But $j x_l^j=x_l\partial/ \partial x_l^j,$
and the result follows immediately. \qed

\begin{lemma}\label{L3}
Let $\alpha$ be a partition, let $l(\alpha)=m+1,$ and let $1\le l\le m.$ Then
$$\sum_{i,j\ge1}x_l^{i+j}\Xi^{\{1,\ldots,m+1\} - \{l,m+1\}}ij \frac{\partial^2}{\partial p_i \partial p_j}  p_\alpha=
\left. x_l\frac{\partial}{\partial x_l}x_{m+1}\frac{\partial}{\partial x_{m+1}}\Xi^{\{1,\ldots,m+1\}}p_\alpha
\right\vert_{x_{m+1}=x_l}.$$
\end{lemma}
\proof Similar to the proof of Lemma~\ref{L2}. \qed

\medskip

The partial differential equation for the symmetrized form of $\Phi$ is given in the following result.

\begin{theorem}\label{T1}
The series $\Psi_m^{(g)},$ for $g,m=1,2,\ldots,$ satisfy the partial differential equation
\begin{eqnarray}\label{pde1}
 \displaystyle{   \left( w_1\frac{\partial}{\partial w_1}+\cdots+w_m\frac{\partial}{\partial w_m}+m+2(g-1)\right)
\Psi_m^{(g)}(x_1,\ldots,x_m) } = T_1+\cdots+T_4,
\end{eqnarray}
where
\begin{eqnarray*}
T_1 &=& \displaystyle{  
\left.\frac{1}{2}\sum_{i=1}^m
\left( x_i\frac{\partial}{\partial x_{i}} x_{m+1}\frac{\partial}{\partial x_{m+1}}\Psi^{(g-1)}_{m+1}(x_1,\ldots,x_{m+1})\right)
\right\vert_{x_{m+1}=x_i} }, \\
T_2 &=&  \displaystyle{  {\sym1} \frac{w_2}{1-w_1}\frac{1}{w_1-w_2} x_1\frac{\partial}{\partial x_{1}}
 \Psi^{(g)}_{m-1}(x_1,x_3,\ldots,x_m)  },\\
T_3 &=& 
\displaystyle{ \sum_{k=3}^m \sym{k-1}\left( x_1\frac{\partial}{\partial x_{1}} \Psi^{(0)}_k(x_1,\ldots,x_k)\right)
  \left(x_1\frac{\partial}{\partial x_{1}}\Psi^{(g)}_{m-k+1}(x_1,x_{k+1},\ldots,x_m)\right) }, \\
T_4 &=&
\displaystyle{  \frac{1}{2}\sum_{\atp{1\le k\le m}{1\le a\le g-1}}\sym{k-1} 
\left(  x_1\frac{\partial}{\partial x_1} \Psi^{(a)}_k(x_1,\ldots,x_k)\right)
  \left( x_1\frac{\partial}{\partial x_1}\Psi_{m-k+1}^{(g-a)}(x_1,x_{k+1},\ldots,x_m)\right)  }.
\end{eqnarray*}
with initial condition~(\ref{Pe1m}), where $\Psi_0^{(g)}=0$ for $g=1,2,\ldots.$
\end{theorem}
\proof Let $E_1$ denote the left hand side of equation~(\ref{eLPhi}), and let $E_2, E_3, E_4$
denote the three summations on the right hand side (indexing them from left to right). We apply
$[z^g]\Xi^{\{1,\ldots,m \}},$ with $u=1,$ to each of these in turn, and then denote the results of this
by $U_2, U_3$ and $U_4,$ respectively. 

\noindent{\bf For $E_1:$} We consider, with $l(\alpha)=m, \alpha\vdash n,$
\begin{eqnarray*}
\Xi^{\{1,\ldots,m \}}(n+m+2g-2) p_\alpha  
&=&\left(\sum_{j=1}^n\alpha_j+m+2g-2\right)\Xi^{\{1,\ldots,m \}}p_\alpha  \\
&=&\left( \sum_{i=1}^mx_i\frac{\partial}{\partial x_i}+m+2g-2\right)\Xi^{\{1,\ldots,m \}}p_\alpha,
\end{eqnarray*}
so
\begin{eqnarray*}
\left.[z^g]\Xi^{\{1,\ldots,m \}} E_1 \right\vert_{u=1} 
&=& \left( x_1\frac{\partial}{\partial x_1}+\cdots+x_m\frac{\partial}{\partial x_m}+m+2g-2\right)
\Psi_m^{(g)}(x_1,\ldots,x_m).
\end{eqnarray*}

\noindent{\bf For $E_2:$} 
We consider, with $l(\alpha)=m+1,$
\begin{eqnarray*}
\Xi^{\{1,\ldots,m \}}\sum_{i,j\ge1}p_{i+j}ij\frac{\partial^2 p_\alpha}{\partial p_i\partial p_j} 
&=& \sum_{l=1}^m \sum_{i,j\ge1} x_l^{i+j} 
\Xi^{\{1,\ldots,m \}-\{l\}}
ij \frac{\partial^2 p_\alpha}{\partial p_i\partial p_j},\\
&=&  \sum_{l=1}^m \sum_{i,j\ge1} x_l^{i+j}
\Xi^{\{1,\ldots,m+1 \}-\{l,m+1\}}
ij \frac{\partial^2 p_\alpha}{\partial p_i\partial p_j},\\
&=&\left. \sum_{l=1}^m \left(
x_l\frac{\partial}{\partial x_l} x_{m+1}\frac{\partial}{\partial x_{m+1}}
\Xi^{\{1,\ldots,m+1 \}}p_\alpha\right)\right\vert_{x_{m+1}=x_l},
\end{eqnarray*}
by Lemma~\ref{L3}. Thus
$$
U_2=\frac{1}{2} \sum_{l=1}^m \left.\left(
x_l\frac{\partial}{\partial x_l} x_{m+1}\frac{\partial}{\partial x_{m+1}}
\Psi^{(g-1)}_{m+1}(x_1,\ldots,x_{m+1}) \right)\right\vert_{x_{m+1}=x_l}.
$$

\noindent{\bf For $E_3:$}
We consider, with $l(\alpha)=k, \,l(\beta)=m-k+1,$
\begin{eqnarray*}
\Xi^{\{1,\ldots,m \}}\sum_{i,j\ge1}p_{i+j}ij\frac{\partial p_\alpha}{\partial p_i}
\frac{\partial p_\beta}{\partial p_j}&=& 
\sum_{l=1}^m \sum_{i,j\ge1} x_l^{i+j}
\Xi^{\{1,\ldots,m \}-\{l\}} \left( i\frac{\partial p_\alpha}{\partial p_i}\right)
\left( j\frac{\partial p_\beta}{\partial p_j}\right) \\
&=& \sum_{l=1}^m
\sum_{\aatp
{ \cA\cup\cB=\{1,\ldots,m \}-\{l\} }
{\vert\cA\vert=k-1,\vert\cB\vert=m-k}
{\cA\cap\cB=\emptyset}
}
\left(\sum_{i\ge1}  x_l^i \Xi^{\cA}\,i\,\frac{\partial p_\alpha}{\partial p_i}\right)\,
\left(\sum_{j\ge1}  x_l^j \Xi^{\cB}\,j\,\frac{\partial p_\beta}{\partial p_j}\right)\\
&\mbox{}&\mbox{by Lemma~\ref{L1}} \\
&=& \sum_{l=1}^m
\sum_{\aatp
{\cA\cup\cB=\{1,\ldots,m \}-\{l\} }
{\vert\cA\vert=k-1,\vert\cB\vert=m-k}
{\cA\cap\cB=\emptyset}
}
\left( x_l\frac{\partial}{\partial x_l} \Xi^{\cA\cup\{l\}}p_\alpha\right)\,
\left( x_l\frac{\partial}{\partial x_l} \Xi^{\cB\cup\{l\}}p_\beta\right)\\
&\mbox{}&\mbox{by Lemma~\ref{L2}} \\
&=&\sym{k-1}
\left( x_1\frac{\partial}{\partial x_1} \Xi^{\{1,\ldots,k \}}p_\alpha\right)\,
\left( x_1\frac{\partial}{\partial x_1} \Xi^{\{1,k+1,\ldots,m \}}p_\beta\right),
\end{eqnarray*}
whence
$$
U_3=\frac{1}{2}\sum_{\atp{1\le k\le m} {0\le a\le g}} \sym{k-1}
\left(x_1\frac{\partial}{\partial x_1} \Psi^{(a)}_k(x_1,\ldots,x_k)\right)\,
\left(x_1\frac{\partial}{\partial x_1} \Psi_{m-k+1}^{(g-a)}(x_1,x_{k+1},\ldots,x_m)\right).
$$

\noindent{\bf For $E_4:$}
We consider, with $l(\alpha)=m-1,$
\begin{eqnarray*}
\Xi^{\{1,\ldots,m \}}\sum_{i,j\ge1}(i+j)p_i p_j\frac{\partial p_\alpha}{\partial p_{i+j}}
 &=& 2\sum_{1\le l<k\le m} \sum_{i,j\ge1}x_l^i x_k^j\,
\Xi^{\{1,\ldots,m \}-\{l,k\}}(i+j)
\frac{\partial p_\alpha}{\partial p_{i+j}} \\
&=& 2\sum_{1\le l<k\le m}\,\, \sum_{r\ge1}
\frac{x_kx_l^r-x_lx_k^r}{x_k-x_l}\,\,
\Xi^{\{1,\ldots,m \}-\{l,k\}}
r\frac{\partial p_\alpha}{\partial p_{r}} \\
 &=& 2 \sum_{1\le l\neq k\le m}  \frac{x_k}{x_l-x_k}\, x_l\frac{\partial}{\partial x_l}
\Xi^{\{1,\ldots,m \}-\{k\}}p_\alpha \quad\mbox{by Lemma~\ref{L2}}\\
&=& 2\sym1  \frac{x_2}{x_1-x_2}\, x_1\frac{\partial}{\partial x_1}
\Xi^{\{1,3,4\ldots,m \}}p_\alpha.
\end{eqnarray*}
Thus
$$
U_4= \sym1 \frac{x_2}{x_1-x_2}\, x_1\frac{\partial}{\partial x_1}
\Psi_{m-1}^{(g)}(x_1,x_3,x_4,\ldots,x_m).$$

Collecting these cases and combining them, we have
\begin{eqnarray}\label{u234}
\left(\sum_{i=1}^m x_i\frac{\partial}{\partial x_i}+m+2(g-1)\right) \Psi_m^{(g)}(x_1,\ldots,x_m)
=U_2+U_3+U_4.
\end{eqnarray}
Now $U_2$ gives $T_1$ in~(\ref{pde1}). Also, from~(\ref{Pe1m}),
 $\left(x_{1}\frac{\partial}{\partial x_{1}}\right)^2\Psi^{(0)}_1=w_1/(1-w_1),$ 
so from~(\ref{yxw}), we obtain $$x_{1}\frac{\partial}{\partial x_{1}}\Psi^{(0)}_1(x_1)=w_1.$$
This allows us to simplify the terms with  $a=0,k=1$ and $a=g,k=m$ in $U_3,$ giving
$$\sum_{i=1}^mw_ix_i \frac{\partial}{\partial x_{i}}\Psi_m^{(g)}(x_1,\ldots,x_m).$$
This summation is then moved to the left hand side of  equation~(\ref{u234}),
and the left hand side of equation~(\ref{pde1}) follows immediately from~(\ref{yxw}).

Also, from~(\ref{Pe1m}), 
$\left(x_{1}\frac{\partial}{\partial x_{1}}+x_{2}\frac{\partial}{\partial x_{2}}\right)\Psi^{(0)}_2=w_1w_2/(1-w_1)(1-w_2),$
and it can be verified that
$\Psi^{(0)}_2(x_1,x_2)=\log((w_1-w_2)/(x_1-x_2))-w_1-w_2,$ whence
\begin{eqnarray}\label{ex1d}
\left(x_{1}\frac{\partial}{\partial x_{1}}\right)\Psi^{(0)}_2(x_1,x_2) =\frac{w_2}{(1-w_1)(w_1-w_2)}-\frac{x_2}{x_1-x_2}.
\end{eqnarray}
This allows us to simplify the terms with $a=0,k=2$ and $a=g,k=m-1$ in $U_3,$ giving
$$\sym1 \left(\frac{w_2}{(1-w_1)(w_1-w_2)}-\frac{x_2}{x_1-x_2}\right)
x_{1}\frac{\partial}{\partial x_{1}}\Psi^{(g)}_{m-1}(x_1,x_3,\ldots,x_m).$$
This expression is then combined with $U_4,$ to give $T_2$ in~(\ref{pde1}).
The terms with $a=0,k=3,\ldots,m$ and $a=g,k=1,\ldots,m-2$ in $U_3$ combine to give
$T_3$ in~(\ref{pde1}). Finally, the remaining terms in $U_3$ give $T_4$ in~(\ref{pde1}), and the
result follows from~(\ref{u234}).  \qed

\section{Determining the symmetric generating series}\label{soi}
We now consider how to use the transformed partial differential equation~(\ref{pde1}) to determine $\Psi_m^{(g)}.$
We  begin with the following result that establishes in general that this series is a rational series in
$w_1,\ldots,w_m.$
\begin{theorem}\label{tPhitd}
For $g\ge1$ and for $g=0, m\ge3,$
$\Psi_m^{(g)}$ is a symmetric polynomial in $(1-w_1)^{-1},\ldots,(1-w_m)^{-1}$
of total degree less than or equal to $2m+6g-5.$
\end{theorem}
\proof The proof is by induction.
The result is true for $g=0,m\ge3,$ from the initial condition~(\ref{Pe1m}),
using~(\ref{wxe}) to express $\Psi^{(0)}_m$ in terms of the $w_i$'s for $i=1,\dots,m$ where
$m\ge3.$

Now for $g=m=1,$ the right hand side of equation~(\ref{pde1})
comes from $T_1$ alone, and is given by
\begin{eqnarray*}
\lim_{x_2\rightarrow x_1}\frac{1}{2} x_1\frac{\partial}{\partial x_1}x_2\frac{\partial}{\partial x_2} \Psi^{(0)}_2(x_1,x_2),
\end{eqnarray*}
which can be evaluated straightforwardly from~(\ref{ex1d}), to give
\begin{eqnarray*}
\frac{1}{24}w_1^2\frac{w_1^2-4w_1+6}{(1-w_1)^4}.
\end{eqnarray*}
Thus for $g=m=1,$ equation~(\ref{pde1}) contains only rational functions of $w_1,$ and the result
is true for $g=m=1.$

To prove the result for $g\ge1, m\ge1,$ $(g,m)\neq(1,1),$ note that equation~(\ref{pde1})
contains no explicit occurrences of the $x_i$'s except through the operator
$x_i\partial / \partial x_i,$ and that~(\ref{yxw}) allows us to replace
these by $(w_i/(1-w_i)) \partial / \partial w_i,$ for $i=1,\ldots,m.$
Moreover, for these values of $g,m,$ the right hand side of~(\ref{pde1})  
does not contain $\Psi^{(0)}_1$ or $\Psi^{(0)}_2.$  Thus if the induction
hypothesis is that the result is true for $\Psi^{(i)}_j,$ where $0\le i \le g,$ and
$0\le j\le m+g-i,$ with $(i,j)\neq (0,1),(0,2),(g,m),$ then equation~(\ref{pde1})
is of the form
\begin{eqnarray}\label{pde2}
\displaystyle{   \left( w_1\frac{\partial}{\partial w_1}+\cdots+w_m\frac{\partial}{\partial w_m}+m+2(g-1)\right)
\Psi^{(g)}_m(x_1,\ldots,x_m) }=
K^{(g)}_m(\bfw),
\end{eqnarray}
where $K^{(g)}_m(\bfw)$ depends only on those $\Psi_j^{(i)}$ to which we have applied
the induction hypothesis. Thus $K^{(g)}_m(\bfw)$ is a rational function
of the $w_i$'s and $w_1-w_2$ perfectly divides by symmetry, so
$K^{(g)}_m(\bfw)$ is a polynomial in  $(1-w_1)^{-1},\ldots,(1-w_m)^{-1}.$
Moreover, from the induction hypothesis, we can easily bound the total degree
of $K^{(g)}_m(\bfw)=T_1+\cdots+T_4.$
First note that $x_i\partial / \partial x_i$ increases the degree of a power of 
$(1-w_i)^{-1}$ by $2.$ Thus, from the induction hypothesis, the terms arising from $T_1$
 have total degree less that or equal to
$2m+6g-5,$ and the terms arising from $T_2,T_3$ and $T_4$ have total degree  less that or equal to
$2m+6g-4.$

But $w_i\partial / \partial w_i$ on the left hand side of~(\ref{pde1}) increases the degree of a power
of $(1-w_i)^{-1}$ by $1.$ Thus the result follows by induction, that $\Psi^{(g)}_m$
has total degree less that or equal to
$2m+6g-5,$ for $g\ge2,$ and for $g=1,m\ge2,$ and the proof is complete. \qed

\medskip

Now we return to the question of actually determining $\Psi^{(g)}_m$ from equation~(\ref{pde1}).
If we determine the $\Psi^{(g)}_m$ in the order specified in the induction above, then equation~(\ref{pde1})
becomes~(\ref{pde2}), where $K^{(g)}_m(\bfw)$ is a known rational function of $\bfw,$
and this partial differential equation is trivial, as we show next.

\begin{theorem}\label{LPhiint}
In the notation of~(\ref{pde2}), 
$$\Psi^{(g)}_m=\int_0^1 K^{(g)}_m(t\bfw)t^{m+2g-3}dt.$$
\end{theorem}
\proof In equation~(\ref{pde2}), substitute $tw_i$ for $w_i$ for $i=1,\ldots,m,$
so $\sum_{i=1}^m w_i \partial / \partial w_i$ becomes $td/dt.$ Now multiply both sides
by $t^{m+2g-3},$ giving
\begin{eqnarray*}
\frac{d}{dt}\left( \left. t^{m+2g-2}\Psi^{(g)}_m\right\vert_{\bfw\mapsto t\bfw} \right)=
t^{m+2g-3}K^{(g)}_m(t\bfw),
\end{eqnarray*}
and the result follows. \qed

\medskip 

In principle, Theorem~\ref{LPhiint} provides an iterative procedure for explicitly determining the 
$\Psi_m^{(g)}.$ We have used Maple in this way to determine $\Psi^{(g)}_m$ for
$g=1,m=1,\ldots,6$ and $g=2,\ldots,5, m=1,\ldots,6-g.$ 
Because of some technical issues that arise in using Maple for this
purpose, we in fact modified the procedure in practice, and sketch the details of this in 
Section~\ref{scomments}. 

In each case where we have made computations  $\Psi_m^{(g)}$ can be written in the form
\begin{eqnarray*}
\Psi^{(g)}_m(\bfx)=
f_m^{(g)}\left(x_1\frac{\partial}{\partial x_1},\ldots,x_m\frac{\partial}{\partial x_m}\right)
V_m(\bfw),
\end{eqnarray*}
where 
$f_m^{(g)}(x_1{\partial}/{\partial x_1},\ldots,x_m{\partial}/{\partial x_m})$
is a symmetric polynomial in 
$x_1{\partial}/{\partial x_1},\ldots,x_m{\partial}/{\partial x_m}.$
We have found that  these polynomials are expressed most compactly
 in terms of the elementary symmetric functions
$e_k\equiv e_k(x_1{\partial}/{\partial x_1},\ldots,x_m{\partial}/{\partial x_m}),$
for $k=1,\ldots,m $ and,
for this reason, the Appendix gives the  values of the $f_m^{(g)},$ that are defined by this process, as polynomials
in $e_1,\ldots,e_m.$
We conjecture next that this is always the case, and that the total degree of the polynomial
is $m+3g-3$ in general, as can be observed in all cases in the Appendix.

\begin{conjecture}\label{cPhi}
For $g=0, m\ge3$ and for $g\ge1,m\ge1,$
\begin{eqnarray}\label{ePhfg}
\Psi^{(g)}_m(\bfx)=
f_m^{(g)}\left(x_1\frac{\partial}{\partial x_1},\ldots,x_m\frac{\partial}{\partial x_m}\right)
V_m(\bfw),
\end{eqnarray}
for a unique symmetric polynomial $f_m^{(g)}(x_1,\ldots,x_m)$  of total degree
$m+3g-3$  in $x_1,\ldots,x_m.$ 
\end{conjecture}

Conjecture~\ref{cPhi} is also in agreement with the previous known results, as given in Section~\ref{ssback}.
For example,
Hurwitz's result~(\ref{er1}) for $g=0,$ as restated in~(\ref{Pe1m}),
 gives $f^{(0)}_m=e_1^{m-3}$ for $m\ge3.$
The form of Conjecture~\ref{cPhi} might seem surprising, compared with the form of the solution established in
Theorem~\ref{tPhitd}. However, as we discuss in Section~\ref{scomments},
Conjecture~\ref{cPhi} is only a slight strengthening of Theorem~\ref{tPhitd}.

Conjecture~\ref{cPhi}  is also  especially useful for determining the ramification numbers
$\mu_m^{(g)}(\alpha)$ for, as we show in the next result, $f_m^{(g)}$
evaluated at the argument $\alpha$ actually gives $\mu_m^{(g)}(\alpha)$
up to a known scaling factor.

\begin{lemma}\label{mucfm}
If~(\ref{ePhfg}) holds, then for $l(\alpha)=m,$ $\alpha=(\alpha_1,\ldots,\alpha_m),$
$$\mu_m^{(g)}(\alpha)=\frac{1}{n!}\vert\cC_\alpha\vert\,(n+m+2g-2)! \left(\prod_{j=1}^m
\frac{\alpha_j^{\alpha_j}}{(\alpha_j-1)!}\right)f_m^{(g)}(\alpha).$$
\end{lemma}  
\proof The solution of~(\ref{wxe}) is
$$w_r=\sum_{j\ge1}\frac{j^{j-1}}{j!}x_r^j,$$
by Lagrange's Theorem, so from~(\ref{Vkpw})
$$\left[x_r^{\alpha_r}\right]\,\left(x_r\frac{\partial}{\partial x_r}\right)^{k_r}
\frac{w_r}{1-w_r} 
=\left[x_r^{\alpha_r}\right]\,\left(x_r\frac{\partial}{\partial x_r}\right)^{k_r+1}w_r
=\frac{\alpha_r^{\alpha_r+k_r}}{\alpha_r!}.$$
Thus if $\alpha=(\alpha_1,\ldots,\alpha_m),$ then
$$\left[\bfx^{\alpha}\right]\,\prod_{r=1}^m\left(x_r\frac{\partial}{\partial x_r}\right)^{k_r}
 V_m(\bfw)=
\prod_{r=1}^m\alpha_r^{k_r}    
\, \frac{\alpha^\alpha}{\alpha!}.$$
so by linearity, we have
\begin{eqnarray}\label{efV}
\left[\bfx^{\alpha}\right]\,f_m^{(g)}\left(x_1\frac{\partial}{\partial x_1},\ldots,x_m\frac{\partial}{\partial x_m}\right)
V_m(\bfw)
=f_m^{(g)}(\alpha) \frac{\alpha^\alpha}{\alpha!}.
\end{eqnarray}

Now we look at $\Psi^{(g)}_m(\bfx)$ as defined in~(\ref{eXi}) and~(\ref{dPhi}).
First, by counting how often each monomial appears, we have
$$\frac{\vert\cC_\alpha\vert}{n!}\sum_{\omega\in\symgp_m}\bfx^{\omega(\alpha)}
=\frac{1}{\prod_{j=1}^m \alpha_j}m_\alpha(\bfx),$$
where $m_\alpha$ is a monomial symmetric function. Thus $\Psi^{(g)}_m(\bfx)$ can be reexpressed in the form
\begin{eqnarray*}
\Psi^{(g)}_m(\bfx)=\sum_{n\ge1}\sum_{\atp{\alpha\vdash n} {l(\alpha)=m}}
\frac{1}{(n+m+2(g-1))!}\frac{1}{\prod_{j=1}^m \alpha_j}c_g(\alpha)m_\alpha(\bfx).
\end{eqnarray*}
Thus
$$\left[\bfx^\alpha\right]\Psi_m^{(g)}(\bfx) = \frac{1}{(n+m+2g-2)!}\prod_{j=1}^m\alpha_j c_g(\alpha),$$
and the result follows by comparing  this with~(\ref{emug}),~(\ref{efV}) and~(\ref{ePhfg}).  \qed

In Section~\ref{Si}, we have recast Conjecture~\ref{cPhi} and Lemma~\ref{mucfm}
to make them refer as immediately as possible to the ramification numbers $\mu^{(g)}_m$:
in Section~\ref{ssback} we have {\em defined} $f^{(g)}_m$ in~(\ref{newequ}) as the rescaling
of $\mu^{(g)}_m$ given in the statement of Lemma~\ref{mucfm};
in Section~\ref{ssback} we give Conjecture~\ref{conject1},
that $f^{(g)}_m$  thus defined is a symmetric polynomial of total degree $m+3g-3.$
Of course, this means that the symmetric polynomials $f^{(g)}_m$ given in the Appendix
thus have two equivalent interpretations.
One interpretation is in the arguments $x_1\partial/\partial x_1,\ldots,x_m\partial/\partial x_m,$
giving a partial differential operator as considered in Conjecture~\ref{cPhi}.
The second interpretation is in the arguments $\alpha_1,\ldots,\alpha_m,$
giving a scaled expression for $\mu^{(g)}_m$ as considered in~(\ref{newequ}).

Finally, as further strengthening of Theorem~\ref{tPhitd}, for $g=1$ we are able to conjecture
a  closed form from the data in the Appendix, as follows.
Note the similarity to Hurwitz's result that $f_m^{(0)}=e_1^{m-3}.$
\begin{conjecture}\label{cfmg}
For  $m\ge1,$
\begin{eqnarray*}
f^{(1)}_m=\frac{1}{24} \left(e_1^m-e_1^{m-1}-\sum_{i=2}^m(i-2)!e_{i}e_1^{m-i}\right).
\end{eqnarray*}
\end{conjecture}

As stated above, the data we have obtained in the Appendix, for $g=1,m=1,\ldots,6$ is in
agreement with Conjecture~\ref{cfmg}. This is restated as Conjecture~\ref{cfm24} because of
its strikingly simple form.

\section{Computational comments}\label{scomments}
We first address the relationship between the forms for $\Psi^{(g)}_m$ proved in Theorem~\ref{tPhitd}
and conjectured in Conjecture~\ref{cPhi}.

Note that, for $j\ge1,$
\begin{eqnarray}\label{extoc}
\left(x_i\frac{\partial}{\partial x_i}\right)^j \frac{w_i}{1-w_i} =
\sum_{k=j+1}^{2j+1}c(j,k)\frac{1}{(1-w_i)^k},
\end{eqnarray}
(for $j=0,$ the lower limit of the summation becomes $j$),
for some integers $c(j,k).$ Also, since
$w_i \partial / \partial w_i= (1-w_i)x_i \partial / \partial x_i,$
we have for $j\ge1,$ from~(\ref{yxw}),
\begin{eqnarray}\label{ewxtoc}
w_i\frac{\partial}{\partial w_i}\left(x_i\frac{\partial}{\partial x_i}\right)^{j-1} \frac{w_i}{1-w_i} &=&
\sum_{k=j}^{2j}c(j,k)\frac{1}{(1-w_i)^k}.
\end{eqnarray}
This triangular system of linear equations can be inverted to express 
$((1-w_i)^{-1})^k$ as a linear combination of 
$$\left(x_i\frac{\partial}{\partial x_i}\right)^j\quad\mbox{and}\quad
w_i\frac{\partial}{\partial w_i}\left(x_i\frac{\partial}{\partial x_i}\right)^j.$$
Consequently, since from~(\ref{Vkpw})
\begin{eqnarray*}
V_m(\bfw)=\prod_{i=1}^m\left( \frac{1}{1-w_i} -1\right),
\end{eqnarray*}
then Theorem~\ref{tPhitd} implies that, in general,
$\Psi_m^{(g)}$ can be written as  $f_m^{(g)}V_m(\bfw),$
where $f_m^{(g)}$ is a differential operator which is
a symmetric sum of monomials in  $x_i\partial/\partial x_i,$ $i=1,\ldots,m$
plus a sum of a single  $w_j\partial/\partial w_j$ multiplied by symmetric sums in
the $x_i\partial/\partial x_i$'s.
Thus Conjecture~\ref{cPhi} strengthens Theorem~\ref{tPhitd} concerning the form of
$\Psi_m^{(g)}$ only in asserting that the terms involving $w_j\partial/\partial w_j,$
for $j=1,\ldots,m$ do not appear. This is in agreement with all known values
for $\Psi_m^{(g)},$
although we
have been unable to prove  that this holds in general. Also in the conjecture,
the degree of this symmetric differential operator follows from~(\ref{extoc})
and the upper bound given in Theorem~\ref{tPhitd}, and it is part of the conjecture that the bound is actually
attained.

In implementing Theorem~\ref{LPhiint} (using {\sf Maple})  to determine explicitly the $\Psi_m^{(g)}$ given
in the Appendix, we have actually used~(\ref{extoc}) and~(\ref{ewxtoc}), and their inverses, to transform the process
into one involving only polynomial  instead of more awkward rational series.
In particular, this process gives a proof of Theorem~\ref{T12}


\section{The proof of Theorem~\ref{T11}}\label{ST11}
We conclude with a proof of Theorem~\ref{T11}. The proof uses the series $w_1$ that has been
introduced for the symmetrization of $\Phi,$ and is the solution of the functional equation~(\ref{wxe}).

According to~\cite{V}, $\mu_n^{(1)}\equiv \mu_n^{(1)}(1^n)$ satisfies the following recurrence
equation that was discovered by Pandharipande and Graber:

$$
\mu_n^{(1)}=\frac n6\binom{n}{2}(2n-1)\mu_n^{(0)}+2(2n-1)\sum_{j=1}^{n-2} (n-j)j^2
\binom{2n-2}{2j-2} \mu_j^{(0)} \mu_{n-j}^{(1)},
$$
where $\mu_n^{(0)}\equiv \mu_n^{(0)}(1^n).$
Since $\mu_n^{(0)}=(2n-2)!n^{n-3}/n!$ and $\mu_n^{(1)}=(2n)!f_n^{(1)}/n!$
from~(\ref{er1}) and~(\ref{newequ}), respectively, where $f_n^{(1)}\equiv f_n^{(1)}(1^n),$
then
$$
\frac{(2n)!}{n!}f^{(1)}_n=\frac n6 \binom{n}{2} \frac{(2n-1)!}{n!} n^{n-3}
+2(2n-1)! \sum_{j=1}^{n-2} \frac{n-j}{j!(n-j)!} j^{j-1} f^{(1)}_{n-j}.
$$
Then $a_n$ defined by $a_n=24n f^{(1)}_n$ satisfies the recurrence equation
$$
a_n = (n-1) n^{n-1} +  \sum_{j=1}^{n-2}\binom{n}{j}  j^{j-1} a_{n-j}.
$$ 

Let $w\equiv w_1.$ Then, by Lagrange's Theorem~(see, for example,~\cite{GJbook})
$$w=\sum_{n\ge1}\frac{n^{n-1}}{n!}x^n \quad\mbox{and}\quad \frac{w}{1-w}=\sum_{n\ge1}\frac{n^{n}}{n!}x^n.$$
Let
$$A(x)=\sum_{n\ge2}a_n\frac{x^n}{n!}.$$
It follows from the recurrence equation for $a_n$ that
$$A(x)=\frac{w}{1-w} - w + wA(x)$$
so
$$A(x)=\frac{w^2}{(1-w)^2}.$$
Now note $A(x)$ can be re-expressed in the form
\begin{eqnarray*}
A(x) &=& \frac{w}{(1-w)^3}-\frac{w}{1-w} -\frac{w^2}{(1-w)^3} \\
&=& \frac{w}{(1-w)^3}-\frac{w}{1-w} -
\sum_{i\ge2}\frac{1}{i(i-1)} \frac{iw^i-(i-1)w^{i+1}} {(1-w)^{3}}.
\end{eqnarray*}
But, again by Lagrange's Theorem,
$$
\frac{w}{(1-w)^3}=\sum_{n\ge1}\frac{n^{n+1}}{n!}x^n\quad\mbox{and}\quad
\frac{iw^i-(i-1)w^{i+1}}{(1-w)^{3}} = \sum_{j\ge0}\frac{(i+j)^{j+1}}{j!}x^{i+j}.
$$
The second of these is obtained by considering the expansion of $w^i/(1-w)$ and
$(xd/dx)(w^i/(1-w)),$ having observed that $w/(1-w)=xdw/dx,$ from the functional equation for $w.$
Each term in the expression for $A(x)$ has now been given a series expansion,
and result follows. 


\section*{Acknowledgements}
This work was supported by grants individually
to IPG and DMJ from the Natural
Sciences and Engineering Research Council of Canada.



\appendix
\section{Explicit results}\label{A1}
The following  gives explicit expressions  for $f_m^{(g)}.$ 
The results can be expressed in a 
compact form in terms of weighted divided differences $\Delta^{(g)}_m$ as follows.

Let $\alpha\vdash n$ and $l(\alpha)=m.$ 
Let $d_1=3!2^2,$ $d_2=6!2^3,$ $d_3=9!2^3,$  and $d_4=12!2^5.$ 
Then, for $g=1,m=1,\ldots,6$;  $g=2,3,4, m=1,\dots,6-g,$

$$d_gf^{(g)}_m=e_1d_gf^{(g)}_{m-1}+e_m\Delta^{(g)}_m,$$
for $m\ge1,$ where $f^{(g)}_0\equiv0,$ $g\ge1.$
Thus
\begin{eqnarray*}
f^{(g)}_m &=& \frac{1}{d_g}
\left(
e_1^{m} \Delta^{(g)}_1 +e_1^{m-2}e_2 \Delta^{(g)}_2 +e_1^{m-3}e_3 \Delta^{(g)}_3 +  \cdots+e_m\Delta^{(g)}_m\right).
\end{eqnarray*}

\begin{eqnarray}\label{egenus1}
\left.
\begin{array}{ccll}
\mbox{{\bf Genus 1:}}&\Delta_1^{(1)}&=&1-\frac{1}{e_1},\\
&\Delta_2^{(1)}&=&-0!,\\
&\Delta_3^{(1)}&=&-1!,\\
&\Delta_4^{(1)}&=&-2!,\\
&\Delta_5^{(1)}&=&-3!,\\
&\Delta_6^{(1)}&=&-4!.
\end{array}\right\}
\end{eqnarray}
\begin{eqnarray}\label{egenus2}
\left.\begin{array}{ccll}
\mbox{{\bf Genus 2:}}&\Delta_1^{(2)}&=& \phantom{-}5e_1^3-12e_1^2+7e_1,\\
&\Delta_2^{(2)}&=&\left(-10e_1^3+9e_1e_{2}\right)+\left(12e_1^2-2e_{2}\right),\\
&\Delta_3^{(2)}&=&\left(-18e_1^3+18e_1e_{2}-3e_{3}\right)+\left(16e_1^2-6e_{2}\right),\\
&\Delta_4^{(2)}&=&\left(-36e_1^3+60e_1e_{2}-12e_{3}\right)+\left(38e_1^2-24e_{2}\right).
\end{array}
\right\}
\end{eqnarray}
\begin{eqnarray}\label{egenus3}
\left. \begin{array}{ccll}
\mbox{{\bf Genus 3:}}&\Delta_1^{(3)}&=&  35e_1^6-147e_1^5+205e_1^4-94e_1^3,\\
&\Delta_2^{(3)}&=& \left(-105 e_1^6 +189e_1^4e_{2} -135e_1^2e_{2}^2\right)
 +\left(294e_1^5 -321e_1^3e_{2} +90e_1e_2^2\right) \\
&&\mbox{}&+\left(-205e_1^4+74e_1^2e_2-16e_2^2\right),\\
&\Delta_3^{(3)}&=&\left(-273e_1^6+594e_1^4e_2 +153e_1^3e_3\right.\\
&&\mbox{}&\left. -405e_1^2e_2^2+135e_1e_2e_3-27e_3^2\right) \\
&&\mbox{}&+ \left(642e_1^5 -912e_1^3e_2-111e_1^2e_3+360e_1e_2^2-66e_2e_3\right)  \\
&&\mbox{}& +\left( -353e_1^4+270e_1^2e_2 +64e_1e_3-80e_2^2 \right).
\end{array}\right\} 
\end{eqnarray}
\begin{eqnarray}\label{egenus4}
\left.\begin{array}{ccll}
\mbox{{\bf Genus 4:}}&\Delta_1^{(4)}&=& 1925e_1^9-12320e_1^8+29854e_1^7-32032e_1^6+12573e_1^5,\\
&\Delta_2^{(4)}&=&\left( -7700e_1^9+20790e_1^7e_2-29700e_1^5e_2^2+17325e_1^3e_2^3\right) \\
&&\mbox{}&+\left(36960e_1^8-74316e_1^6e_2+72600e_1^4e_2^2-23100e_1^2e_2^3\right)\\
&&\mbox{}&+\left(-59708e_1^7 +77814e_1^5e_2-44880e_1^3e_2^2+10780e_1e_2^3\right)\\
&&\mbox{}& +\left(32032e_1^6-182e_1^4e_2+8800e_1^2e_2^2-1584e_2^3\right).
\end{array}\right\}
\end{eqnarray}



\begin{thebibliography}{999}
\bibitem{A1}{\sc V.Arnold}
{\em  Critical points of functions and classification of caustics},
Uspekhi Mat. Nauk.,
{\bf29} (1994), 243--244.

\bibitem{Arntcc}{\sc V.Arnold}
{\em  Topological classification of complex trigonometric polynomials and the combinatorics of   graphs
with the same number of edges and vertices},
Functional Anal, Appl.,
{\bf30} (1996), 1--17, (Russian).

\bibitem{CT}{\sc M.Crescimanno and W.Taylor},
{\em  Large N phases of chiral $QCD_2$},
Nuclear Phys. B,
{\bf437} (1995), 3--24.

\bibitem{ELSV}{\sc T.Ekedahl, S.Lando, M.Shapiro, and A.Vainshtein},
{\em  Hurwitz numbers and characteristic classes of the moduli spaces},
(1998), (in preparation).

\bibitem{EK}{\sc L.Ernstrom and G.Kennedy},
{\em Recursive formulas for characteristic numbers of rational plane curves},
J.Algebraic Geom.,
{\bf7} (1998), 141--181.

\bibitem{Gdosf}{\sc I.P.Goulden},
{\em  A differential operator  for symmetric functions and the combinatorics
of multiplying  transpositions},
Trans. Amer. Math. Soc.,
{\bf344} (1994), 421--440.

\bibitem{GJtfst}{\sc I.P.Goulden and D.M.Jackson},
{\em Transitive factorizations in the symmetric group and combinatorial aspects of singularity theory},
(July 1997), (preprint).

\bibitem{GJtransf}{\sc I.P.Goulden and D.M.Jackson},
{\em Transitive factorizations into transpositions and holomorphic mappings on the sphere},
Proc. Amer. Math. Soc.,
{\bf125} (1997), 51--60.

\bibitem{GJsfmr}{\sc  I.P.Goulden and D.M.Jackson},
{\em Symmetric functions and Macdonald's result for top connection coefficients in the symmetric group},
J.~Algebra,
{\bf166} (1994), 364--378.

\bibitem{GJbook}{\sc I.P.Goulden and D.M.Jackson}
``Combinatorial Enumeration,''
Wiley, New York,
1983.

\bibitem{GL}{\sc V.Goryunov and S.Lando},
{\em  On enumeration of meromorphic functions on the line},
 (1997), (preprint).

\bibitem{Hrfgv}{\sc A.Hurwitz},
{\em Ueber Riemann'sche Fl\"{a}chen mit gegebenen Verzweigungspunkten},
Matematische Annalen,
{\bf39} (1891), 1--60.

\bibitem{dmjcppcc}{\sc D.M.Jackson},
{\em  Some combinatorial problems associated with products of conjugacy classes
of the symmetric group},
J.~Combinatorial Theory,
{\bf49} (1988), 363--369.

\bibitem{KZ}{\sc A.Khovanskii and S.Zdravkovska},
{\em Branched covers of $\sS^2$ and braid groups},
J.~Knot Theory and Its Ramifications,
{\bf5} (1996), 55--75.

\bibitem{K}{\sc M.Kontsevich},
{\em  Enumeration of rational curves via torus actions},
in ``The Moduli Space of Curves,'' R.Dijkgraaf a.o. eds., Birkhauser, 
(1995), pp.~335--368.

\bibitem{KM}{\sc M.Kontsevich and Yu.Manin},
{\em  Gromov-Witten classes, quantum cohomology, and enumerative geometry},
Comm. Math. Phys.,
{\bf164} (1994), 525--562.

\bibitem{L}{\sc E.Looijenga},
{\em The complement of the bifurcation variety of a simple singularity},
Invent. Math.,
{\bf23} (1973), 105--116.

\bibitem{MAC}{\sc I.G.Macdonald},
``Symmetric functions and Hall polynomials,''
Second Edition, Clarendon Press, Oxford, 1995.

\bibitem{medn}{\sc A.D.Mednykh},
{\em  Branched coverings of Riemann surfaces whose branch orders coincide
with the multiplicity},
Comm. in Algebra.,
{\bf18} (1990), 1517--1533.

\bibitem{N}{\sc S.Natanzon},
{\em  Topology of 2-dimensional coverings and meromorphic functions on real and complex algebraic curves},
Selecta Mathematica (formerly Sovetica),
{\bf12} (1993), 251--291.

\bibitem{P}{\sc R.Pandharipande},
{\em  Intersection of $\rats$-divisors on Kontsevich's moduli space $\overline{M}_{0,n}(\sP^r,d)$
and enumerative geometry},
Trans. Amer. Math. Soc.,
(1998), (to appear).


\bibitem{ssv}{\sc B.Shapiro, M.Shapiro and A.Vainshtein},
{\em  Ramified coverings of $S^2$ with one degenerate branching point and
enumeration of edge-ordered graphs},
Adv. in Math. Sci., {\bf34} (1997), 219--228.

\bibitem{vs}{\sc V.Strehl},
{\em  Reconstruction of a proof by A.Hurwitz},
(November 1996), (private communication).

\bibitem{V}{\sc R.Vakil},
{\em  Enumerative geometry of plane curves of low genus},
AG-9803007.

\bibitem{Z}{\sc D.Zvonkine},
{\em  Multiplicities of the Lyashko-Looijenga map and its strata},
C.R. Acad. Sci., Paris,
{\bf324} (1997), 1349--1353.


\end{thebibliography}
\end{document}